\documentclass[12pt]{article}
\usepackage{graphicx}
\usepackage{amsmath,amsthm,amssymb,enumerate}
\usepackage{euscript,mathrsfs}
\usepackage{color}
\usepackage{dsfont}
\usepackage[left=2cm,right=2cm,top=3.5cm,bottom=3.5cm]{geometry}
\usepackage{color}
\usepackage[framemethod=tikz]{mdframed}
\usepackage{bm}
\allowdisplaybreaks

\usepackage{soul}

\catcode`\@=11 \@addtoreset{equation}{section}

\catcode`\@=12

\allowdisplaybreaks

\newtheorem{Theorem}{Theorem}[section]
\newtheorem{Proposition}[Theorem]{Proposition}
\newtheorem{Lemma}[Theorem]{Lemma}
\newtheorem{Corollary}[Theorem]{Corollary}

\theoremstyle{definition}
\newtheorem{Definition}[Theorem]{Definition}

\newtheorem{Remark}[Theorem]{Remark}

\newcommand{\bTheorem}[1]{
\begin{Theorem} \label{T#1} }
\newcommand{\eT}{\end{Theorem}}

\newcommand{\bProposition}[1]{
\begin{Proposition} \label{P#1}}
\newcommand{\eP}{\end{Proposition}}

\newcommand{\bLemma}[1]{
\begin{Lemma} \label{L#1} }
\newcommand{\eL}{\end{Lemma}}

\newcommand{\bCorollary}[1]{
\begin{Corollary} \label{C#1} }
\newcommand{\eC}{\end{Corollary}}

\newcommand{\bRemark}[1]{
\begin{Remark} \label{R#1} }
\newcommand{\eR}{\end{Remark}}

\newcommand{\bDefinition}[1]{
\begin{Definition} \label{D#1} }
\newcommand{\eD}{\end{Definition}}

\newcommand{\rhN}{\vr^{h,N}}
\newcommand{\uhN}{\vu^{h,N} }

\newcommand{\tmu}{\widetilde{\mu}}
\newcommand{\teta}{\widetilde{\eta}}

\newcommand{\vrh}{\vr^h}

\newcommand{\Td}{\mathbb{T}^d}
\newcommand{\bFormula}[1]{
\begin{equation} \label{#1}}
\newcommand{\eF}{\end{equation}}

\newcommand{\vuh}{\vu^h}

\newcommand{\Ov}[1]{\overline{#1}}

\newcommand{\aleq}{\stackrel{<}{\sim}}

\newcommand{\vr}{\varrho}

\newcommand{\tvr}{\tilde \vr}
\newcommand{\tvu}{{\tilde \vu}}

\newcommand{\vu}{\vc{u}}
\newcommand{\vm}{\vc{m}}

\newcommand{\vc}[1]{{\bm #1}}

\newcommand{\Div}{{\rm div}_x}
\newcommand{\Grad}{\nabla_x}

\newcommand{\dx}{\,{\rm d} {x}}

\newcommand{\dt}{\,{\rm d} t }

\newcommand{\intTd}[1]{\int_{\mathbb{T}^d} #1 \ \dx}

\newcommand{\D}{{\rm d}}

\newcommand{\ep}{\varepsilon}

\newcommand{\expe}[1]{ \mathbb{E} \left[ #1 \right] }

\newcommand{\br}{ \nonumber \\ }

\def\softd{{\leavevmode\setbox1=\hbox{d}%
          \hbox to 1.05\wd1{d\kern-0.4ex{\char039}\hss}}}
\definecolor{Cgrey}{rgb}{0.85,0.85,0.85}
\definecolor{Cblue}{rgb}{0.50,0.85,0.85}
\definecolor{Cred}{rgb}{1,0,0}
\definecolor{fancy}{rgb}{0.10,0.85,0.10}

\newcommand\Cbox[2]{%
    \newbox\contentbox%
    \newbox\bkgdbox%
    \setbox\contentbox\hbox to \hsize{%
        \vtop{
            \kern\columnsep
            \hbox to \hsize{%
                \kern\columnsep%
                \advance\hsize by -2\columnsep%
                \setlength{\textwidth}{\hsize}%
                \vbox{
                    \parskip=\baselineskip
                    \parindent=0bp
                    #2
                }%
                \kern\columnsep%
            }%
            \kern\columnsep%
        }%
    }%
    \setbox\bkgdbox\vbox{
        \color{#1}
        \hrule width  \wd\contentbox %
               height \ht\contentbox %
               depth  \dp\contentbox
        \color{black}
    }%
    \wd\bkgdbox=0bp%
    \vbox{\hbox to \hsize{\box\bkgdbox\box\contentbox}}%
    \vskip\baselineskip%
}

\mdfdefinestyle{MyFrame}{%
	linecolor=black,
	outerlinewidth=1pt,
	roundcorner=5pt,
	innertopmargin=\baselineskip,
	innerbottommargin=\baselineskip,
	innerrightmargin=10pt,
	innerleftmargin=10pt,
	backgroundcolor=white!20!white}

\date{}

\makeindex
\begin{document}


\title{Compressible fluid motion with uncertain data}

\author{Eduard Feireisl
	\thanks{The work of E.F. was partially supported by the
		Czech Sciences Foundation (GA\v CR), Grant Agreement
		21--02411S. The Institute of Mathematics of the Academy of Sciences of
		the Czech Republic is supported by RVO:67985840. }
}

\date{\today}

\maketitle

\bigskip

\centerline{$^*$  Institute of Mathematics of the Academy of Sciences of the Czech Republic}

\centerline{\v Zitn\' a 25, CZ-115 67 Praha 1, Czech Republic}

\begin{abstract}
	
	We propose a suitable analytical framework to perform numerical analysis of problems arising in 
	compressible fluid models with uncertain data. We discuss both weak and strong stochastic approach, where the former is based on the knowledge of the mere distribution (law)  of the random data typical for the Monte-Carlo and related methods, while the latter assumes the data to be known as a random variable 
	on a given probability space aiming at obtaining 
	the associated solution in the same form. As an example of the strong approach, we discuss the stochastic 
	collocation method based on a piecewise constant approximation of the random data.

\end{abstract}

\bigskip

\centerline{\it Dedicated to the memory of my friend Anton\' \i n Novotn\' y}


\section{Introduction}
\label{i}

Standard mathematical models in continuum fluid mechanics based on the Navier--Stokes system suffer the well 
known deficiencies confirmed recently by the results of Buckmaster and Vicol \cite{BucVic1}, Buckmaster, 
Cao--{L}abora, G{\' o}mez--{S}errano
\cite{BuCLGS}, Merle et al. \cite{MeRaRoSz}:
\begin{itemize}
\item Weak solutions exist globally in time but may not be uniquely determined by the data.	
\item Strong solutions are uniquely determined by the data but their life--span may be finite.
	\end{itemize} 

We focus on problems involving compressible and (linearly) viscous fluids governed by the 
barotropic Navier--Stokes system. The abstract framework, however, is applicable to a larger class of equations, in particular those including thermal effects.
Despite the above mentioned rather pessimistic scenarios, problems related to viscous fluid are globally well posed at least for smooth initial data close to an equilibrium solutions, cf. Matsumura and Nishida \cite{MANI}, 
Valli and Zajaczkowski \cite{VAZA}, among others. In addition, we anticipate  that smoothness is a generic property 
for a statistically significant set of data. Accordingly, we propose a theoretical framework to perform numerical analysis of problems with uncertain data and ``generic'' well posedness. Our principal working 
hypothesis is the \emph{boundedness in probability} of the approximate numerical solutions specified 
and discussed in Section \ref{B}. Very roughly indeed one may see it as a statistical counterpart of convergence of bounded numerical solutions to deterministic problems. In the present context, the crucial role 
in the analysis plays the conditional regularity property of strong solutions to the compressible Navier--Stokes system 
proved by Sun, Wang, and Zhang \cite{SuWaZh}. Note that boundedness of approximate solutions is considered to be a very mild and acceptable hypothesis used explicitly in a number of numerical studies. 

\subsection{Model problem}

The motion of a compressible (linearly) viscous fluid is described in terms of its mass density $\vr = \vr(t,x)$ and the 
velocity field $\vu = \vu(t,x)$ satisfying the following equations.

\begin{mdframed}[style=MyFrame]
	
	\textsc{Navier--Stokes system.}
	
	\begin{align}
		\partial_t \vr + \Div (\vr \vu) &= 0, \label{i1}\\
		\partial_t (\vr \vu) + \Div (\vr \vu \otimes \vu) + \Grad p &= \Div \mathbb{S} (\Grad \vu) + 
		\vr \vc{g},
		\label{i2} \\
		\mathbb{S}(\Grad \vu) &= \mu \left( \Grad \vu + \Grad^t \vu - \frac{2}{d} \Div \vu \mathbb{I} \right) +
		\eta \Div \vu \mathbb{I},\ \mu > 0,\ \eta \geq 0. \label{i3}
		\end{align}
	
	\end{mdframed}

\noindent
We refer e.g. to Gallavotti \cite{GALL1} for the physical background of the model.

Eliminating the effect of thermal changes we close the system by imposing the barotropic (isentropic) equation of state for the pressure, 
\begin{equation} \label{i4}
p = p(\vr) = a \vr^\gamma,\ a > 0,\ \gamma > 1.	
	\end{equation}

For the sake of simplicity, we identify the spatial domain with the flat torus corresponding 
to the space periodic boundary conditions, while the motion is considered on a given compact time interval.

\begin{mdframed}[style=MyFrame]
	
	\textsc{Periodic boundary conditions.}

\begin{equation} \label{i5}
\ t \in [0,T],\	x \in \Td, \ d=2,3.
	\end{equation}

\end{mdframed}

Finally, the initial state of the system is prescribed.

\begin{mdframed}[style=MyFrame]
	
\textsc{Initial data.}	
\begin{equation} \label{i6}
	\vr(0, \cdot) = \vr_0,\ \vr \vu (0, \cdot) = \vr_0 \vu_0,
	\end{equation}

\end{mdframed}

As a matter of fact, the periodic boundary conditions are considered only for simplicity, similar results can be obtained, for instance, for the physically relevant no--slip condition 
\[
\vu|_{\partial Q} = 0,
\]
where $Q$ represent the physical domain occupied by the fluid. Note, however, that many numerical methods are 
applied on polygonal (non--smooth) domains, where the problem does not admit smooth solutions required in the subsequent analysis. The numerical method must be therefore adapted to approximate not only the exact solutions but also the physical space $Q$, see e.g. \cite{FeHoMaNo2015_I}, \cite{FeHoMaNo}. 

\subsection{Data}

The solution $(\vr, \vu)$ of the problem \eqref{i1}--\eqref{i6} is determined by the data:

\begin{itemize}
	
	\item the initial data $\vr_0$, $\vu_0$;
	\item the viscosity coefficients $\mu$ and $\eta$;
	\item the parameter $a$ in the equation of state;
	\item the driving force $\vc{g}$.
	
	\end{itemize}

The set of data 
\[
D = \left[ \vr_0, \vu_0, \mu, \eta, a, \vc{g} \right] \in X_D
\]
will be a subset of the data space $X_D$. Ideally, $X_D$ is a separable Hilbert space, in particular Polish, 
to make  the standard tools of stochastic analysis applicable.
The initial data as well as the driving force will be regular as we require the problem to admit a smooth (classical) solution at least locally in time.

One is tempted to say that the subsequent analysis may apply to more general data including, in particular, different types of boundary conditions studied by Kwon and Novotn\' y \cite{KwoNovII}, \cite{KwoNovI}.  Unfortunately, the blow-up criteria based on boundedness of strong solutions are not available so far for general boundary value problems.

\subsection{Data uncertainty, statistical solution}

Our aim is to study the problem \eqref{i1}--\eqref{i6} with random (uncertain) data. Accordingly, the 
mapping 
\[
t \in [0,T] \times [\vr_0, \vu_0, \mu, \eta, a, \vc{g} ] \to (\vr, \vu)(t,\cdot) 
\]
can be considered as a stochastic process solving the system of equations a.s. Such a process is 
usually termed \emph{statistical solution}. We refer to Constantin and Wu \cite{ConWu}, Foias et al. \cite{FoMoTi}, \cite{FoRoTe1}, Vishik and Fursikov \cite{VisFur} for the relevant theory in the context of incompressible fluids. The corresponding ``compressible'' analogue is presented in \cite{FanFei}. Note, however, that these results are based on the concept of weak solutions existing globally in time.

Our goal is to study convergence of suitable numerical approximations in both weak and strong statistical settings. In the \emph{weak} setting, only statistical properties of the data as well as of the associated approximate and exact solutions are relevant. These are, for instance, the mean (expected value), the median, r-means of higher order etc. 

The \emph{strong} setting is based on the exact representation of the data as a random variable, with the goal to obtain the solution in the same form. Needless to say that all properties available in the weak setting may be recovered knowing the exact form of the strong solution. Here, we illustrate this approach by analysing the collocation method presented in \cite{FeiLuk2021}.

The paper is organized as follows. In Section \ref{E}, we recall the necessary ingredients of the mathematical theory of the compressible Navier--Stokes system. In Section \ref{N}, we introduce the concept of numerical approximation and the associated numerical methods. We do not specify the method contenting ourselves with stating its basic properties required by the present theory. In Section \ref{R}, we discuss in detail the random data related to the weak and strong stochastic approach. In Section \ref{B}, we introduce our main working hypothesis concerning the boundedness in probability of the approximate solutions. Finally, in Section \ref{C}, we establish convergence of the approximate numerical solutions and state the main results.

\section{The existence theory for the Navier--Stokes system}
\label{E}

We recall the well known facts concerning the existence and regularity of strong solutions to the compressible Navier--Stokes system.

\begin{mdframed}[style=MyFrame]
	
	\begin{Theorem}[{\bf Local existence}] \label{ET1}
	Let $k \geq 3$.
		Suppose the data belong to the class 
		\begin{align} \label{E1}
			\vr_0 \in W^{k,2}(\Td), \ \inf \vr_0 = \underline{\vr} > 0,\  & 
			 \vu_0 \in W^{k,2}(\Td; R^d), \\
			\label{E2}
			0 < \underline{\mu}&\leq \mu,\ \eta \geq 0,\ 0 < \underline{a} \leq a \leq \Ov{a}, \\
			\vc{g} &\in W^{k + 1,2}((0,\infty) \times \Td; R^d). \label{E3} 
			\end{align}

Then there exists $T_{\rm max} > 0$ such that the Navier--Stokes system \eqref{i4}--\eqref{i6} admits 
a classical solution $\vr$, $\vu$, unique in the class 
\begin{equation} \label{E4}
\vr \in C([0,T]; W^{k,2}(\Td)),\ \vu \in C([0,T]; W^{k,2}(\Td; R^d)) \cap L^2(0,T; W^{k+1,2}(\Td; R^d)),\ 
0 < T < T_{\rm max}.	
	\end{equation} 		
		
		\end{Theorem}

	\end{mdframed}

If $k=3$, Theorem \ref{ET1} was proved by Valli and Zajaczkowski \cite{VAZA}. Regularity for $k \geq 3$ was established in \cite{BrFeHo2016}, cf. also Gallagher \cite{Gall2000} or Tani \cite{TAN}. The proof is based on transforming the Navier--Stokes system to a parabolic perturbation of a symmetric hyperbolic system therefore 
using the specific form of the isentropic pressure equation of state. Moreover, as shown in \cite{BrFeHo2016}, 
\begin{equation} \label{E5}
	T_{\rm max} < \infty \ \Rightarrow \ \limsup_{t \to T_{\rm max}} \| \vu(t,\cdot) \|_{W^{2, \infty}(\Td; R^d)}
	\to \infty.
\end{equation}

In addition, using the remarkable regularity criterion of Sun, Wang, and Zhang \cite{SuWaZh1}, 
the authors in \cite[Proposition 2.2]{FeHoMaNo2015_I} observed that there exists a function 
\[
 C \left( \| \vr_0 \|_{C^3}, \| \vu_0 \|_{C^3}, \underline{\vr}^{-1}, \underline{\mu}, 
\underline{a}^{-1}, \Ov{a},\ \| \vc{g} \|_{C^1},\ 
\sup_{t \in [0, \tau]} \| \vr(t, \cdot) \|_{L^\infty(\Td)}, \sup_{t \in [0, \tau]} \| \vu(t, \cdot), \|_{L^\infty(\Td; R^d)}\right)
\]
which is bounded for bounded values of its arguments, such that
\begin{equation} \label{E6}
\sup_{t \in [0,\tau] } \| \vu(t,\cdot) \|_{W^{2, \infty}(\Td; R^d)} \leq C 
\ \mbox{for all}\  \tau \in [0,T_{\rm max}).
\end{equation}
In particular, the strong solution exists as long as we control its norm in $L^\infty$. 
Thus we have obtained the following conditional regularity result. 
\begin{mdframed}[style=MyFrame]
	
	\begin{Theorem}[{\bf Regularity criterion}] \label{ET2}
		
		Let $k \geq 5$. Let $(\vr, \vu)$ be the local solution of the Navier--Stokes system claimed 
	in Theorem \ref{ET1}. 
	
	Then 
	\begin{align} 
	&\sup_{t \in [0, T]} \left( \| \vr (t, \cdot) \|_{W^{k,2}(\Td)} + \| \vu (t, \cdot) \|_{W^{k,2}(\Td; R^d)} \right) + 
	\int_0^T \| \vu \|^2_{W^{k+1}(\Td; R^d)} \dt \br
	\leq 
	 &C \Big( T, \| (\vr_0, \vu_0) \|_{W^{k,2}}, \underline{\vr}^{-1}, \underline{\mu}, 
	\underline{a}^{-1}, \Ov{a},\ \| \vc{g} \|_{W^{k+1,2}((0,T) \times \Td; R^d)},\ 
	\| (\vr, \vu) \|_{L^\infty((0,T) \times \Td; R^{d + 1})}\Big)
			\label{E7}	
		\end{align}
for any $0 < T < T_{\rm max}$, where $C$ is a bounded function of bounded arguments. In particular, 
\begin{equation} \label{E8}
	T_{\rm max} < \infty \ \Rightarrow \ \limsup_{t \to T_{\rm max}} \| (\vr, \vu)(t,\cdot) \|_{L^\infty (\Td; R^{d+1})}
\to \infty.	
	\end{equation}		
		
		\end{Theorem}
	
\end{mdframed}

\begin{Remark} \label{ER1}
	
	Boundedness of the velocity can be omitted in \eqref{E7}, \eqref{E8} provided $\eta = 0$, cf. 
	\cite{SuWaZh}.
	
	\end{Remark}

Theorem \ref{ET2} is a remarkable result and plays the key role in the subsequent analysis. Indeed a smooth 
solution exists as long as we control its $L^\infty-$norm in perfect agreement with the recently obtained blow--up results \cite{BuCLGS}, \cite{MeRaRoSz}.

In the light of the above arguments, a suitable \emph{data space} $X_D$ is 
\begin{align} 	
\left[ \vr_0, \vu_0, \mu, \eta, a, \vc{g} \right] 	
\in X_D = W^{k,2}(\Td) \times W^{k,2}(\Td; R^d) \times R \times R \times R \times W^{k+1,2}((0,T) \times \Td; R^d),\label{E9}
\end{align}
with $k \geq 5$. 
Clearly, $X_D$ is a separable Hilbert space. In addition, we consider a closed convex subset $A_D \subset X_D$ of 
\emph{admissible data}, 
\begin{equation} \label{R1}
	A_D = \left\{ \left[\vr_0, \vu_0, \mu, \eta, a, \vc{g} \right] \in X_D \ \Big| \ 
	\vr_0 \geq \underline{\vr} > 0,\ \mu \geq \underline{\mu} > 0,\ 
	\eta \geq 0, \ 0 < \underline{a} \leq a \leq \Ov{a},\ \| \vc{g} \|_{L^\infty} \leq \Ov{g}
	\right\}.
\end{equation}

\begin{Remark} \label{RR1}
	The reader will have noticed that our hypotheses are taylored to the physically relevant case $d=3$ and could be possibly relaxed if $d=1,2$.
	\end{Remark}

\section{Numerical approximation}
\label{N}

We are ready to introduce the concept of \emph{numerical approximation} associated to given data 
\begin{align} 
\left[ \vr_0, \vu_0, \mu, \eta, a, \vc{g} \right] &\in A_D, \br
0 < \underline{\vr} &\leq \vr_0, \ \mu \geq \underline{\mu} > 0,\ 0 < \underline{a} \leq a \leq \Ov{a}, 
\label{N1}
\end{align}
namely $\vrh$, $\vuh$, $h = h(\ell) \searrow 0$ as $\ell \to \infty$. 

We suppose that the 
numerical approximation enjoys the following properties:

\begin{itemize}
\item {\bf Admissibility.} $(\vr^h, \vu^h) \in V_h$, where $V_h \subset L^\infty((0,T) \times \Td); R^{d+1})$ is a finite dimensional 
space, 
\[
\inf \vr^h > 0 \ \mbox{for any}\ h, 
\]
\begin{equation} \label{N1a}
\mathcal{A} \left(h, \left[ \vr_0, \vu_0, \mu, \eta, a, \vc{g} \right], \vr_h, \vu_h \right) = 0, 	
	\end{equation}
where 
\[
\mathcal{A} : (0,\infty) \times A_D \times V_h \to R^m ,\ m = m(h) 
\]
is a Borel measurable (typically continuous) mapping representing a finite system of algebraic equations called \emph{numerical scheme}. Note that the scheme \eqref{N1a} may admit several solutions for given data.

\item {\bf Bounded graph property.}

If $N = N(\ell) \nearrow \infty$, $h = h(\ell) \searrow 0$, 
\[
\left[ \vr^N_0, \vu^N_0, \mu^N, \eta^N, a^N, \vc{g}^N \right] \in A_D 
\to \left[ \vr_0, \vu_0, \mu, \eta, a , \vc{g} \right] \ \mbox{in}
\ X_D \ \mbox{as}\ N \to \infty, 
\]
and the associated numerical approximation satisfies
\begin{equation} \label{N2}
\sup_{h,N} \left\| (\rhN, \uhN ) \right\|_{L^\infty((0,T) \times \Td; R^{d+1})} < \infty,
	\end{equation}	
then 
\[
\rhN \to \vr \ \mbox{in}\ L^1((0,T) \times \Td),\ 
\uhN \to \vu \ \mbox{in}\ L^1((0,T) \times \Td; \Td) \ \mbox{as}\ h \to 0,\ N \to \infty,
\]
where $(\vr, \vu)$ is the unique classical solution of the Navier--Stokes system associated to the data 
$\left[ \vr_0, \vu_0, \mu, \eta, a, \vc{g} \right]$.

	\end{itemize}

Very roughly indeed, the bounded graph property requires any bounded numerical approximation to be convergent. 
Apparently, this depends on the character of the numerical scheme represented by the operator $\mathcal{A}$ in 
\eqref{N1a}. Validity of the bounded graph property for a time implicit finite volume scheme is established  
in \cite{FeiLuk2021}, cf. also the monograph \cite[Theorem 11.3]{FeLMMiSh}. 

It is interesting to note that any \emph{convergent} numerical approximation of the Navier--Stokes system 
actually admits the open graph property. 

\begin{Definition}[{\bf Convergent numerical approximation}] \label{ND1}
	
We say that a numerical approximation of the Navier--Stokes system \eqref{i1}--\eqref{i6} is \emph{convergent}
if for any sequence of data 
\[
[\vr_0^N, \vu^N_0, \mu^N, \eta^N, \vc{g}^N] \in A_D \to 
[\vr_0, \vu_0, \mu, \eta, \vc{g} ] \ \mbox{in}\ X_D \ \mbox{as}\ N \to \infty,
\]
the numerical apporximation $(\vr^{h,N}, \vu^{h,N})$ satisfies: 
\begin{itemize}\item
	\[
	\vr^{h,N} > 0; 
	\]
	
	\item 
	\[
	\vr^{h,N} \to \vr \ \mbox{in}\ L^1((0,T) \times \Td),\ 
	\vu^{h,N} \to \vu \ \mbox{in}\ L^1((0,T) \times \Td; R^d) \ \mbox{as}\ N \to \infty, \ h \to 0, 
	\]
	for any $0 < T < T_{\rm max}$, where $(\vr, \vu)$ is the unique classical solution of the problem 
	with the data $[\vr_0, \vu_0, \mu, \eta, \vc{g}]$ defined on the maximal time interval $[0, T_{\rm max})$.

	\end{itemize}
	
	\end{Definition}

As a direct consequence of the regularity criterion established in Theorem \ref{ET2}, we get the following

\begin{mdframed}[style=MyFrame]

\begin{Corollary} \label{NT1}
	
	Any convergent numerical approximation is admissible \\ and possesses the bounded graph property.
	
	\end{Corollary}

\end{mdframed}

\section{Random data}

\label{R}

In any real implementation of a numerical scheme, the random data must be approximated by a finite number of deterministic ones. Recall that any data belong to the set 
\[
	A_D = \left\{ \left[\vr_0, \vu_0, \mu, \eta, a, \vc{g} \right] \in X_D \ \Big| \ 
	\vr_0 \geq \underline{\vr} > 0,\ \mu \geq \underline{\mu} > 0,\ 
	\eta \geq 0, \ 0 < \underline{a} \leq a \leq \Ov{a},\ \| \vc{g} \|_{L^\infty} \leq \Ov{g}
	\right\},
	\]
where the constants $\underline{\vr}$, $\underline{\mu}$, $\underline{a}$, $\Ov{a}$, $\Ov{g}$ are \emph{deterministic}.

\subsection{Weak stochastic approach}
\label{WS}

Weak stochastic approach, typical for the Monte Carlo and related methods, is based on the knowledge 
of \emph{random data distribution} - their law in the Polish space $X_D$. Given a random variable $D$ ranging in a Polish space $X_D$, we denote its distribution (law) 
$\mathcal{L}[D]$ -- a Borel probability measure on the space $X_D$. 
The expected output of the method is formulated analogously in terms of the distribution of the associated 
exact solution or, more precisely, its numerical approximation. 

In practice, it is possible to generate 
statistical \emph{samples} of the random data with only finitely many items. Given the (random) data 
$[\vr_0, \vu_0, \mu, \eta, a, \vc{g}] \in A_D$, our starting point is generating a sequence 
of (deterministic) approximations
\begin{equation} \label{R2}
 \left[ \vr_{0}^{n}, \vu_{0}^{n}, \mu^{n}, \eta^{n}, a^n, \vc{g}^n\right] \in A_D	
	\end{equation}
such that for any $F \in BC(X_D)$ there holds
\begin{equation} \label{R3}
\frac{1}{N}	\sum_{n=1}^N F \left[\vr_{0}^{n}, \vu_{0}^{n}, \mu^{n}, \eta^{n}, a^n, \vc{g}^n \right] \to 
\expe{ F [\vr_0, \vu_0, \mu, \eta, a, \vc{g}] } \ \mbox{as}\ N \to \infty, 
	\end{equation}
where the expected value on the right--hand side is given as
\[
\expe{ F [\vr_0, \vu_0, \mu, \eta, a, \vc{g}] } =
\int_{X_D} F \left( \hat{\vr}, \hat{\vu}, \hat{\mu} , \hat{\eta}, \hat{a}, \hat{\vc{g}} \right) \ \D \mathcal{L}[\vr_0, \vu_0, \mu, \eta, a, \vc{g}]. 
\]
Associating to $\left[\vr_{0,n}, \vu_{0,n} \mu_n, \eta_n, a_n, \vc{g}_n \right]$ the sequence 
of discrete (empirical) measures 
\[
\mathcal{L}_N = \frac{1}{N} \sum_{n=1}^N \delta_{\left[\vr_{0}^{n}, \vu_{0}^{n}, \mu^{n}, \eta^{n}, a^n, \vc{g}^n \right]} \in \mathfrak{P}(X_D),
\]
we may equivalently reformulate \eqref{R3} as 
\begin{equation} \label{R4}
\mathcal{L}_N \to \mathcal{L}{[\vr_0, \vu_0, \mu, \eta, a, \vc{g}]} 
\ \mbox{weakly in}\ \mathfrak{P} (X_D) \ 
 \mbox{as}\ N \to \infty
\end{equation}
where $\mathfrak{P}$ denotes the set of Borel probability measures on $X_D$.

Leaving apart the problem of specific \emph{construction} of a suitable approximate sequence \\ $\left[\vr^{0,n}, \vu^{0,n}, \mu^n, \eta^n, a^n, \vc{g}^n \right]$ satisfying \eqref{R3}, we consider its
numerical approximation \\ $(\vr^{h,n}, \vu^{h,n})$ specified in Section \ref{N}.

The goal is to establish a convergence result for the sequence of measures
\[
\frac{1}{N} \sum_{n=1}^N \delta_{(\vr^{h,n}, \vu^{h,n})}.
\]
The numerical solutions $(\vr^{h,n}, \vu^{h,n} )$ belong {\it a priori} to the finite dimensional space 
$V_h \subset L^\infty$, however, $L^\infty$ is not a separable space. Therefore it is more convenient to consider the negative spaces $W^{-m,2}((0,T) \times \Td) \times W^{-m,2}((0,T) \times \Td ; R^d)$ that are separable Hilbert spaces. The desired result would then be 
\begin{equation} \label{R6}
	\frac{1}{N} \sum_{n=1}^N {F} [ \vr^{h,n} , \vu^{h,n} ] \to \expe{ F[\vr, \vu] } \ \mbox{as}\ h \to 0, \ N \to \infty,	
\end{equation}
for any $F \in BC \Big( W^{-m,2}((0,T) \times \Td) \times W^{-m,2}((0,T) \times \Td; R^d) \Big)$, where 
$m > d + 1$, and $(\vr, \vu)$ is the classical solution of the Navier--Stokes system corresponding to the 
data $[\vr_0, \vu_0, \mu, \eta, a, \vc{g}]$. Similarly to the above, the right--hand side of \eqref{R6} 
is interpreted as 
\[
\expe{ F[\vr, \vu] } = \int_{X_D} F \left[ (\vr, \vu) [\hat{\vr}, \hat{\vu}, \hat{\mu} , \hat{\eta}, \hat{a}, \hat{\vc{g}}  ] \right]    \D \mathcal{L}[\vr_0, \vu_0, \mu,  \eta, a, \vc{g} ], 
\]
where $(\vr, \vu)[\hat{\vr}, \hat{\vu}, \hat{\mu} , \hat{\eta}, \hat{a}, \hat{\vc{g}} ]$ is the exact solution associated to the data $[\hat{\vr}, \hat{\vu}, \hat{\mu} , \hat{\eta}, \hat{a}, \hat{\vc{g}}]$. 
Of course, a rigorous justification of \eqref{R6} requires existence of the exact solution on the whole time 
interval $(0,T)$ $\mathcal{L}[\vr_0, \vu_0, \mu, \eta, a , \vc{g}]$- a.s.

If more estimates on the moments of the numerical approximation are available, we may anticipate the convergence 
of empirical means, specifically
\begin{equation} \label{R5}
\frac{1}{N} \sum_{n=1}^N ( \vr^{h,n} , \vu^{h,n} ) \to \expe{ \vr, \vu } \ \mbox{as}\  N \to \infty,\ h \to 0	
	\end{equation}
in a suitable topology, say, $L^q((0,T) \times \Td; R^{d + 1})$, $q \geq 1$. The limit is understood as the Bochner intergral
\[
\expe{ \vr, \vu } = \int_{X_D}  (\vr, \vu) [ \hat{\vr}, \hat{\vu}, \hat{\mu} , \hat{\eta}, \hat{a}, \hat{\vc{g}}   ] \  \D \mathcal{L}[\vr_0, \vu_0, \mu, \lambda, \eta, s, \vc{g} ]
\]
in a suitable Banach space. The relevant rigorous results are stated in Theorems \ref{CT1}, \ref{CT5} in Section \ref{C} below.

\begin{Remark} \label{RR2}
	
	Neither the approximate sequence $[\vr_{0}^n, \vu^n_0, \mu^n, \eta^n, a^n, \vc{g}^n]$ nor the associated numerical solutions $(\vr^{h,n}, \vu^{h,n})$ are  uniquely determined by the data $[\vr_0, \vu_0, \mu, \eta, a, \vc{g}]$. As a matter of fact, the practical implementation deals with a large number 
	of \emph{samples} -- sequences $[\vr_{0}^n, \vu^n_0, \mu^n, \eta^n, a^n, \vc{g}^n]$ -- generated independently 
	mimicking the Strong law of large numbers, cf. e.g. Mishra, Schwab et al.  
	\cite{Koley}, \cite{Kuo}, \cite{Leonardi}, \cite{Mishra_Schwab}.
	
	\end{Remark}

\subsection{Strong stochastic approach}
\label{SS}

Strong stochastic approach requires the precise knowledge of the data as a random variable 
\[
[\vr_0, \vu_0, \mu, \eta, a, \vc{g}] : \left\{ \Omega, \mathcal{B}, \mathcal{P} \right\} 
\to X_D.
\]
defined on a probability space $\Omega$, with a family of measurable sets $\mathcal{B}$, and a complete 
probability measure $\mathcal{P}$. 
The goal is to identify an approximate numerical solution as a random variable \emph{in the same probability space}.

There many ways of suitable approximation, here
we focus on the statistical \emph{collocation method} proposed in \cite{FeiLuk2021}. Writing 
\[
\Omega = \cup_{n=1}^N \Omega_n^N,\ \Omega^N_n \ \mathcal{P}-\mbox{measurable}, \ \Omega^N_i \cap \Omega^N_j = \emptyset \ \mbox{for}\ i \ne j,\ \cup_{n = 1}^N \Omega^n_N = \Omega,
\]
we consider the data
\[
\left[\vr_{0,N}, \vu_{0,N}, \mu_{N}, \eta_{N}, a_N, \vc{g}_N \right] = 
\sum_{n=1}^N \mathds{1}_{\Omega^n_N} (\omega) [\vr_0, \vu_0, \mu, \eta, a, \vc{g}] (\omega_n),\ 
\omega_n \in \Omega^n_N.
\]
Instead of the weak convergence \eqref{R3}, we require the strong convergence of the data,
\begin{equation} \label{R7}
\sum_{n=1}^N \mathds{1}_{\Omega^n_N} (\omega) [\vr_0, \vu_0, \mu, \eta, a, \vc{g}] (\omega_n) \to 
[\vr_{0}, \vu_{0}, \mu, \eta, a, \vc{g}] \ \mbox{in}\ X_D \ \mathcal{P}-\mbox{a.s.}
\end{equation}
We point out that validity of \eqref{R7} may depend on the partition $(\Omega^N_n)_{n=1}^N$ as well as on the position of the collocation points $\omega_n$.
Sufficient conditions for \eqref{R7} to hold can be found in \cite{FeiLuk2021}.

Similarly to the preceding section, we associate to each set of data $[\vr_0, \vu_0, \mu, \eta, a, \vc{g}] (\omega_n)$ its numerical approximation 
\[
[\vr_0, \vu_0, \mu, \eta, a, \vc{g}] (\omega_n) \mapsto [\vr^{h,n}, \vu^{h,n}], 
\]
and a sequence of random variables 
\[
\sum_{n=1}^N \mathds{1}_{\Omega^n_N}(\omega) (\vr^{h,n}, \vu^{h,n})(t,x).
\]
Our goal is to establish the convergence 
\begin{equation} \label{R8}
	\sum_{n=1}^N 1_{\Omega_N^n} (\vr^{h,n}, \vu^{h,n}) \to (\vr, \vu) \ 
	\mbox{as}\ N \to \infty,\ h \to 0\ \mathcal{P}-\mbox{a.s.}
		\end{equation}
in a suitable topology,	
where $(\vr, \vu)$ is the classical solution of the Navier--Stokes system \eqref{i1}--\eqref{i6} corresponding to the data $[\vr_0, \vu_0, \mu, \eta, \vc{g}]$. The relevant results are stated in Theorems \ref{CT2}, \ref{CT4}
below.

\section{Boundedness in probability of approximate solutions}
\label{B}

The crucial \emph{hypothesis} we impose on the family of approximate numerical solutions is its \emph{boundedness in 
	probability}, which can be seen as a discrete counterpart of the bounded graph property introduced in Section \ref{N}. 
Keeping in mind the constructive aspect of numerical approximation we restrict to countably many approximate solutions. 
More specifically, we suppose
\[
h = h(\ell),\ N = N(\ell),\ h(\ell) \searrow 0,\ N(\ell) \nearrow \infty \ \mbox{as}\ \ell \to \infty.
\]	
This sequential structure will be assumed here and hereafter.

\subsection{Weak stochastic approach}

Given the initial data as in Section \ref{WS}, the approximate solution 
\[
\frac{1}{N} \sum_{n=1}^N \delta_{[\vr^{h,n}, \vu^{h,n}]}
\] 
can be identified with a discrete probability measure on the space $V_h$. In this context, the hypothesis of 
boundedness of probability can be formulated as follows.

\begin{mdframed}[style=MyFrame]
	
	\textsc{Boundedness in probability (weak).}
	
	For any $\ep > 0$, there is $M = M(\ep)$ such that 
	\begin{equation} \label{B1}
		\frac{ \# \left\{ \| \vr^{h,n}, \vu^{h,n} \|_{L^\infty((0,T) \times \Td; R^{d +1}) } > M, \ n\leq N \right\} }{N} < \ep\ 
		\mbox{for any}\ \ell = 1,2,\dots
		\end{equation}

	\end{mdframed}

\subsection{Strong stochastic approach}

In terms of the strong stochastic approach introduced in Section \ref{SS}, the analogue of 
\eqref{B1} reads 

\begin{mdframed}[style=MyFrame]
	
	\textsc{Boundedness in probability (strong).}
	
	For any $\ep > 0$, there is $M = M(\ep)$ such that 
	\begin{equation} \label{B2}
		\sum_{ n \leq N, \left\{ \| \vr^{n,h}, \vu^{n,h} \|_{L^\infty((0,T) \times \Td; R^{d +1}) } > M \right\} } |\Omega_n^N| < \ep
		\ \mbox{for}\ \ell = 1,2,\dots
	\end{equation}

\end{mdframed}

Apparently, condition \eqref{B2} is a ``weighted'' version of \eqref{B1}, with the weights proportional to 
the expectation of the discrete events. Below we show that the convergence of the numerical solutions in the weak and strong approach can be studied in a universal setting by means of a convenient representation of 
the weak random data. 

\section{Convergence of numerical approximations}
\label{C}

We start by recalling several tools from the theory of probability.

\subsection{Basic tools of the theory of probability}

\begin{Theorem} [{\bf Skorokhod's representation theorem}] \label{CTT1}
	
	Let $(\nu_N )_{N = 1}^\infty$ be a sequence of Borel probability measures on a Polish space $X$ such that 
	\[
	\nu_N \to \nu \ \mbox{weakly (narrowly) in}\ \mathfrak{P}(X), 
	\]
	meaning
	\[
	\int_X F (y) \D \nu_N \to \int_X F(y) \D \nu \ \mbox{for any}\ F \in BC(X).
	\]
	
	Then there is a probability space $\{ \Omega; \mathcal{B}, \mathcal{P} \}$ and a sequence of random 
	variables 
	\[
	Y_N : \Omega \to X,\ \mathcal{L}[Y_N] = \mu_N
	\]
	such that 
	\[
	Y_N \to Y \ \mbox{in}\ X \ \mathcal{P}-\mbox{a.s.},\ \mathcal{L}[Y] = \nu.
	\]
	
	\end{Theorem}

\begin{Theorem}[{\bf Prokhorov's theorem}] \label{CTT2}
	
	Let $(\nu_N)_{N=1}^\infty$ be a family of probability measures on a Polish space $X$. 
	
	The following is equivalent:
	
	\begin{itemize}
		
		\item $(\nu_N)_{N=1}^\infty$ is weakly precompact, meaning there is a subsequence such 
		\[
		\nu_{N_k} \to \nu \ \mbox{weakly in}\ \mathfrak{P}(X).
		\]
		
		\item $(\nu_N)_{N=1}^\infty$ is tight, meaning for any $\ep > 0$, there is a compact set $K(\ep) \subset X$ such that 
		\[
		\mu_N (K) \geq 1 - \ep \ \mbox{for all}\ N=1,2,\dots.
		\]

		\end{itemize}
	
	\end{Theorem}

\begin{Theorem}[{\bf Gy\"ongy--Krylov convergence criterion}] \label{CTT3}
	
	Let $X$ be a Polish space and $(Y_N)_{N=1}^\infty$ a sequence of $X-$valued random variables.
	
	Then $(Y_N )_{N = 1}^\infty$ converges in probability if and only if for any sequence of joint laws of
	\[
	(Y_{M_k}, Y_{N_k})_{k = 1}^\infty
	\]
	there exists a further subsequence that converge weakly to a probability measure $\nu$ on $X \times X$ such that
	\[
	\nu \left[ (x,y) \in X \times X,\ x = y \right] = 1.
	\]

	\end{Theorem}

\subsection{From weak to strong}

We transform the weakly converging sequence of data distribution identified in Section \ref{WS} to a strongly converging sequence of random variables on a suitable probability space. 
This is a direct consequence of Skorokhod's representation theorem. Indeed there exists a probability space 
$\{ \Omega, \mathcal{B}, \mathcal{P} \}$ and a sequence of random data 
$[\vr_{0,N}, \vu_{0, N}, \mu_N, \eta_N, a_N, \vc{g}_N ] \in A_D$ such that 

\begin{align} 
	\mathcal{L}	[\vr_{0,N}, \vu_{0,N}, \mu_N, \eta_N, a_N, \vc{g}_N ]  &=  \frac{1}{N} \sum_{n=1}^N \delta_{\left[\vr_{0}^{n}, \vu_{0}^{n}, \mu^{n}, \eta^{n}, a^n, \vc{g}^n \right]}, \br
	[\vr_{0,N}, \vu_{0,N}, \mu_N, \eta_N, a_N, \vc{g}_N ] &\to [\tvr_{0}, \tvu_{0}, \widetilde{\mu}, \widetilde{\eta}, \widetilde{a}, \widetilde{\vc{g}} ] \ \mbox{in}\ X_D \ \D y - \mbox{a.s.},
			\label{C2}
\end{align}	
where
\begin{equation} \label{C3}
	{[\tvr_{0}, \tvu_{0}, \widetilde{\mu}, \widetilde{\eta}, \widetilde{a}, \widetilde{\vc{g}} ]}
\sim {[\vr_{0}, \vu_{0}, {\mu}, {\eta}, {a}, {\vc{g}} ]},
\end{equation}
where $\sim$ denotes equality in law.
The convergence of Skorokhod's representation of weak data is therefore the same as in \eqref{R7}. 

In the convergence proof below, we therefore suppose we are given a sequence of (finitely distributed) random data 
\begin{align}
[\vr_{0,N}, \vu_{0,N}, \mu_N, \eta_N, a_N, \vc{g}_N ] &\in A_D, \br 
[\vr_{0,N}, \vu_{0,N}, \mu_N, \eta_N, a_N, \vc{g}_N ] &\to 
[\vr_{0}, \vu_{0}, \mu, \eta , a, \vc{g} ] \ \mbox{in}\ X_D \ \mathcal{P}-\mbox{a.s.}
\label{C5}
\end{align}
along with the associated family of numerical approximations $(\rhN, \uhN)_{N=1}^\infty \in V_h$.
Moreover, in agreement with the hypotheses \eqref{B1}, \eqref{B2}, we suppose the numerical solutions are 
bounded in probability:

For any $\ep > 0$, there is $M(\ep)$ such that 
\begin{equation} \label{C4}
	\mathcal{P} \left\{ \left\| \rhN, \uhN \right\|_{L^\infty((0,T) \times \Td; R^{d+1} } \geq M    \right\} \leq \ep.
	\end{equation}

\subsection{Pointwise boundedness -- another application of Skorokhod's theorem}

We apply once more Skorokhod's theorem to transform boundedness in probability to boundedness a.s. To this end, 
we consider a new sequence that consists of the data, together with the associated numerical solutions augmented by the 
$L^\infty$-norm of the latter. Specifically, 
\[
Y_{h,N} = \Big\{ [\vr_{0,N}, \vu_{0,N}, \mu_N, \lambda_N, a_N, \vc{g}_N ]; (\rhN, \uhN) ; 
\Lambda_{h,N} \Big\},\ \ \mbox{with}\ \Lambda_{h,N} = \| \rhN, \uhN  \|_{L^\infty},
\]
considered as a sequence of random variables ranging in the Polish space
\[
X = X_D \times W^{-m,2}((0,T) \times \Td; R^{d + 1}) \times R, \ m > d + 1.
\]

In accordance with the hypotheses \eqref{C5}, \eqref{C4}, and the compact embedding \\
$L^\infty \hookrightarrow W^{-m,2}((0,T) \times \Td )$, the family of associated (joint) laws 
$\mathcal{L}[Y_{h,N}]$ is tight in $X$. Thus applying Prokhorov's compactness criterion and Skorokhod's 
representation theorem we deduce there is a new probability space $\{ \widetilde{\Omega}; \widetilde{\mathcal{B}}; \widetilde{\mathcal{P}} \}$ and a (sub)sequence of random variables 
\[
\widetilde{Y}_{h_k, N_k} = \left\{ [\tvr_{0,N_k}, \tvu_{0,N_k}, \tmu_{N_k}, \teta_{N_k}, \widetilde{a}_{N_k}, \widetilde{\vc{g}}_N ]; \left( \tvr^{h_k,N_k}, \tvu^{h_k,N_k} \right); \widetilde{\Lambda}_{h_k, N_k} \right\} \in X
\]
such that 
\begin{align}
\Big\{ & [\tvr_{0,N_k}, \tvu_{0,N_k}, \tmu_{N_k}, \teta_{N_k}, \widetilde{a}_{N_k}, \widetilde{\vc{g}}_N ]; \left( \tvr^{h_k,N_k}, \tvu^{h_k,N_k} \right); \widetilde{\Lambda}_{h_k, N_k}\Big\} \br 
&\sim \Big\{ [\vr_{0,N_k}, \vu_{0,N_k}, \mu_{N_k}, \eta_{N_k}, {a}_{N_k}, {\vc{g}}_{N_k} ]; \left( \vr^{h_k,N_k}, \vu^{h_k,N_k} \right), 	
\Lambda_{h_k,N_k} \Big\},\br 
[\tvr_{0,N_k}, \tvu_{0,N_k}, \tmu_{N_k}, \teta_{N_k}, \widetilde{a}_{N_k}, \widetilde{\vc{g}}_{N_k} ] &\to 
[\tvr_0, \tvu_0, \tmu, \teta, \widetilde{a}, \widetilde{\vc{g}}] 
\ \mbox{in} \ X_D \ \widetilde{\mathcal{P}} -\mbox{a.s.}, \br
\mbox{where}\ [\tvr_0, \tvu_0, \tmu, \teta, \widetilde{a}, \widetilde{\vc{g}}] &\sim [\vr_0, \vu_0, \mu, \eta, {a}, {\vc{g}}] \br
\left( \tvr^{h_k, N_k}, \tvu^{h_k,N_k} \right) &\to \left( \tvr, \tvu \right) \ \mbox{in}\ W^{-m,2}((0,T) \times \Td; R^{d+1}) \ \widetilde{\mathcal{P}} -\mbox{a.s.},
\label{C6a}
\end{align}
and
\begin{equation} \label{C6}
\widetilde{\Lambda}_{h_k,N_k} = \| (\tvr^{h_k,N_k}, \tvu^{h_k,N_k}) \|_{L^\infty} \to \widetilde{\Lambda} \ 
\widetilde{\mathcal{P}} -\mbox{a.s.}
	\end{equation}

As the numerical method admits the open graph property introduced in Section \ref{N}, we conclude 
 \begin{equation} \label{C8}
\left( \tvr_{h_k, N_k}, \tvu^{h_k,N_k} \right) \to (\tvr, \tvu) \ \mbox{strongly in}\ L^q ((0,T) \times \Td; R^{d+1}) \ 
\ \widetilde{\mathcal{P}} -\mbox{a.s.} \ \mbox{for any}\ 1 \leq q < \infty,
	\end{equation}
where $(\tvr, \tvu)$ is the unique (statistical) solution of the Navier--Stokes system \eqref{i1}--\eqref{i6}
associated to the data $[\tvr_0, \tvu_0, \tmu, \teta, \widetilde{a}, \widetilde{\vc{g}}]$.
In particular, the Navier--Stokes system \eqref{i1}--\eqref{i6} admits global classical solution of 
the random data $[\tvr_0, \tvu_0, \tmu, \teta, \widetilde{a}, \widetilde{\vc{g}}]$ $\mathcal{P}-$a.s.

Finally, as the limit is uniquely determined by the data, there is no need of subsequences, and we may use the Gy\" ongy--Krylov convergence criterion to conclude that convergence takes place in the original 
probability space:
\begin{equation} \label{C9}
	\left( \rhN, \uhN \right) \to (\varrho, \vu ) \ \mbox{in}\ L^q((0,T) \times \Td; R^{d+1}) 
	\ \mbox{in}\ \mathcal{P}-\mbox{probability},
\end{equation}	
where $(\vr, \vu)$ is the unique solution of the Navier--Stokes system \eqref{i1}--\eqref{i6} corresponding to the data 
$[\vr_0, \vu_0, \mu, \lambda, \vc{g} ]$.

\subsection{Convergence -- summary}

Summing up the results obtained in the preceding section we are ready to formulate the main results 
concerning weakly and strongly converging data.

\begin{mdframed}[style=MyFrame]

\begin{Theorem}[{\bf Convergence of numerical approximations -- weak data}] \label{CT1}
	
	Fix 
	\[
	N = N(\ell) \nearrow \infty,\ h = h(\ell) \searrow 0,\ \ell \to 0,
	\]
	Let 
	\[
	[\vr_{0}^n, \vu_{0}^n, \mu^n, \eta^n, \vc{g}^n ] \in A_D,\ n=1,2,\dots 
	\]
be a sequence of deterministic data satisfying
\[
\frac{1}{N} \sum_{n=1}^N F [\vr_{0}^n, \vu_{0}^n, \mu^n, \eta^n, \vc{g}^n ] \to 
\int_{X_D} F (\hat{\vr}, \hat{\vu}, \hat{\mu} , \hat{\eta}, \hat{a}, \hat{\vc{g}}) \ \D \mathcal{L}_0 \ \mbox{for any}\ F \in BC(X_D), 
\]
for some Borel probability measure $\mathcal{L}_0$ on $X_D$. Let $(\vr^{h,n}, \vu^{h,n})$ be the associated numerical approximations resulting from a convergent numerical method in the sense of Definition \ref{ND1}. 
Suppose that the family of numerical solutions is bounded in probability, meaning 
for  any $\ep > 0$, there is $M = M(\ep)$ such that 
\[
	\frac{ \# \left\{ \| \vr^{h,n}, \vu^{h,n} \|_{L^\infty((0,T) \times \Td; R^{d +1}) } > M, \ n\leq N \right\} }{N} < \ep\ 
	\mbox{for any}\ \ell = 1,2,\dots 
\]	

Then there exists a probability space $\{ \Omega, \mathcal{B}, \mathcal{P} \}$ 
and random data $[\vr_0, \vu_0, \mu_0, \eta_0, \vc{g}_0 ] \in A_D$ such that the following holds:
\begin{itemize}
	\item 
	The Navier--Stokes system \eqref{i1}--\eqref{i6} admits a classical solution $(\vr, \vu)$ on the time interval 
	$(0,T)$ corresponding to the data $[\vr_0, \vu_0, \mu, \eta, a, \vc{g}]$ $\mathcal{P}$--a.s.
	\item 
	\[
	\mathcal{L} [\vr_0, \vu_0, \mu, \eta, a, \vc{g}] = \mathcal{L}_0.
	\]
	\item 
	\begin{equation} \label{C10}
	\frac{1}{N} \sum_{n=1}^N F \left( \vr^{h,n}, \vu^{h,n} \right) \to 
	\int_{X_D} F \Big( (\vr , \vu ) [\hat{\vr}, \hat{\vu}, \hat{\mu} , \hat{\eta}, \hat{a}, \hat{\vc{g}} ] \Big) \D \mathcal{L}_0\ 
	\mbox{as} \ \ell \to \infty
		\end{equation}
for any $F \in BC (L^q((0,T) \times \Omega; R^{d + 1} ) )$,  $1 \leq q < \infty$.	
	
	\end{itemize}

	\end{Theorem}

\end{mdframed}

As for the strongly converging data, we claim the following result.

\begin{mdframed}[style=MyFrame]
	
	\begin{Theorem}[{\bf Convergence of numerical approximations -- strong data}] \label{CT2}
		Fix
		\[
		N = N(\ell) \nearrow \infty,\ h = h(\ell) \searrow 0,\ \ell \to 0.
		\]
		Let  $[\vr_0, \vu_0, \mu, \eta, a, \vc{g}] \in A_D$ be a random variable defined on a probability space 
		$\{ \Omega, \mathcal{B}, \mathcal{P} \}$. Suppose there is a sequence of deterministic data 
		\[
		[\vr_{0}^n, \vu_{0}^n, \mu^n, \eta^n, \vc{g}^n ] \in A_D,\ n=1,2,\dots 
		\]
		and a family of decompositions of $\Omega$, 
		\[
		\Omega = \cup_{n=1}^N \Omega_n^N,\ \Omega^N_n \ \mathcal{P}-\mbox{measurable}, \ \Omega^N_i \cap \Omega^N_j = \emptyset \ \mbox{for}\ i \ne j,\ \cup_{n = 1}^N \Omega^n_N = \Omega,
		\]
		such that 
		\[
		\sum_{n=1}^N \mathds{1}_{\Omega^N_n} [\vr_{0}^n, \vu_{0}^n, \mu^n, \eta^n, \vc{g}^n ] \to 
		[\vr_0, \vu_0, \mu, \eta, a, \vc{g}] \ \mathcal{P}-\mbox{a.s.} 
		\]
		Let $(\vr^{h,n}, \vu^{h,n})$ be the associated numerical approximations resulting from a convergent numerical method in the sense of Definition \ref{ND1}. 
		Suppose that the family of numerical solutions is bounded in probability, meaning 
		for  any $\ep > 0$, there is $M = M(\ep)$ such that 
		\[
			\sum_{ n \leq N, \left\{ \| \vr^{n,h}, \vu^{n,h} \|_{L^\infty((0,T) \times \Td; R^{d +1}) } > M \right\} } |\Omega_n^N| < \ep
		\ \mbox{for any}\ \ell  = 1,2,\dots
		\]

		Then the Navier--Stokes system \eqref{i1}--\eqref{i6} admits a classical solution $(\vr, \vu)$ on the time interval 
			$(0,T)$ corresponding to the data $[\vr_0, \vu_0, \mu, \eta, a, \vc{g}]$ $\mathcal{P}$--a.s., 
			and 
		\begin{equation} \label{C11}
		\sum_{n=1}^N \mathds{1}_{\Omega^N_n} (\vr^{h,n} \vu^{h,n} ) \to (\vr, \vu) 
		\ \mbox{in}\ L^q((0,T) \times \Td; R^{d + 1}) \ \mbox{as}\ \ell \to \infty,  \ 1 \leq q < \infty 
		\ \mbox{in probability.}	
			\end{equation}

	\end{Theorem}

\end{mdframed}

\subsection{Boundedness of higher order moments}
\label{H}

The results stated in Theorem \ref{CT1}, \ref{CT2} are optimal with respect to the hypotheses imposed on the 
data: Convergence in law in the weak setting and convergence a.s. in the strong setting. Better results are expected provided higher moments of the data and/or numerical approximations are bounded. Here, the main issue is that higher integrability imposed on the data may not lead to similar integrability properties 
of the approximate/exact solutions. This problem is related to the {\it a priori} bounds available for 
solutions to the Navier--Stokes system. 

The energy balance for the Navier--Stokes system reads 
\begin{align} 
	\intTd{ \left[ \frac{1}{2} \vr |\vu|^2 + P(\vr) \right](\tau, \cdot) } &+ 
	\int_0^\tau \intTd{ \mathbb{S}(\mu, \eta, \Grad \vu) : \Grad \vu } \dt \\
	&=\intTd{ \left[ \frac{1}{2} \vr_0 |\vu_0|^2 + P(\vr_0) \right] } + \int_0^\tau \intTd{ \vr \vc{g} \cdot \vu } \dt, \label{H1}
	\end{align} 
where $P'(\vr) \vr - P(\vr) = p(\vr)$. Using H\" olders inequality and Gronwall's argument we check easily the estimate
\begin{align} 
\| \vr(	\tau, \cdot) \|^{\gamma}_{L^\gamma(\Td)} + \| \vr \vu (\tau, \cdot) \|^{\frac{2 \gamma}{\gamma + 1} }_{L^{\frac{2 \gamma}{\gamma + 1}}(\Td; R^d) } &\aleq 
	\intTd{ \left[ \frac{1}{2} \vr |\vu|^2 + P(\vr) \right](\tau, \cdot) } \br
	&\leq C(T, \Ov{g} )\left( 1 + \intTd{ \left[ \frac{1}{2} \vr_0 |\vu_0|^2 + P(\vr_0) \right] } \right)
	\label{H2}
	\end{align}
If the approximate numerical solutions satisfy a variant of the inequality \eqref{H2}, then boundedness of moments of the initial energy transfer directly to moments of the numerical approximations at any time $\tau$.
For instance, the finite volume scheme studied in \cite{FeiLuk2021} enjoys this property. 

As a corollary, we get the following improvement of convergence of strong approximations. 

\begin{mdframed}[style = MyFrame]
	
	\begin{Theorem}[{\bf Strong convergence in expectations}] \label{CT4}
		
		In addition to the hypotheses of Theorem \ref{CT2}, suppose that the energy of the numerical solutions 
		is bounded in expectations, meaning 
		\[
		\sum_{n=1}^N |\Omega^M_n| \intTd{\left[ \frac{1}{2} \vr^{h,n} |\vu^{h,n}|^2 + P(\vr^{h,n}) \right] (\tau,  \cdot)    } \aleq 1 \ \mbox{for}\ \tau \in (0,T),\ \ell = 1,2,\dots
		\]
		
	Then 
	
	\[
	\expe{ \left\| \sum_{n=1}^N \mathds{1}_{\Omega^N_n} \vr^{h,n} - \vr \right\|_{L^\gamma((0,T) \times \Td)}^r }
	\to 0 \ \mbox{as}\ \ell \to \infty \ \mbox{for any}\ 1 \leq r < \gamma,
	\]
	\[
	\expe{ \left\| \sum_{n=1}^N \mathds{1}_{\Omega^N_n} \vr^{h,n} \vu^{h,n} - \vr \vu \right\|_{L^{\frac{2 \gamma}{\gamma + 1}} ((0,T) \times \Td; R^d)}^s }
	\to 0 \ \mbox{as}\ \ell \to \infty \ \mbox{for any}\ 1 \leq s < \frac{2 \gamma}{\gamma + 1}.
	\]
	\end{Theorem}

	\end{mdframed}

In order to formulate an analogous result for the weakly converging data, we introduce the concept of \emph{$r$-barycenter} $\mathbb{E}_r[Y]$ of  
a random variable $Y$ defined on a Polish space $(X; d_X)$: 
\[
\mathbb{E}_r[Y] \in X,\ 
\expe{ d_X \left( Y; \mathbb{E}_r[Y] \right)^r } = \min_{Z \in X} \expe{ d_X \left( Y; Z \right)^r },\ r \geq 1,
\]
meaning 
\[
E_r(Y) = {\rm arg} \min_{Z \in X} \expe{ d_X \left( Y; Z \right)^r }.
\]
If $X = L^q((0,T) \times \Td; R^d)$ and $1 < r < \infty$, then 
\begin{itemize}
	\item there exists a unique $r-$barycenter for any $Y$, $\expe{ \| Y \|_{L^q}^r } < \infty$,
	\item $\mathbb{E}_ r[Y]$ depends only on the distribution (law) of $Y$,
		\end{itemize}
see Cuesta and Matran \cite{CueMat}. Thus adapting \cite[Theorem 4]{CueMat} to the present setting, we get the 
following result.	

\begin{mdframed}[style = MyFrame]
	
	\begin{Theorem}[{\bf Strong convergence of barycenters}] \label{CT5}
		
		In addition to the hypotheses of Theorem \ref{CT1}, suppose that the energy of the numerical solutions 
		is bounded in expectations, meaning 
		\[
		\frac{1}{N} \sum_{n=1}^N  \intTd{\left[ \frac{1}{2} \vr^{h,n} |\vu^{h,n}|^2 + P(\vr^{h,n}) \right] (\tau,  \cdot)    } \aleq 1 \ \mbox{for}\ \tau \in (0,T),\ \ell = 1,2,\dots
		\]
		
		Then 
		
		\begin{itemize}
			\item 
			\[
			\frac{1}{N} \sum_{n=1}^N \vr^{h,n} \to \expe{ \vr } \ \mbox{in} \ L^\gamma ((0,T) \times \Td), 
			\]
		\[
		\frac{1}{N} \sum_{n=1}^N \vr^{h,n} \vu^{h,n} \to \expe{ \vr \vu } \ \mbox{in} \ L^{\frac{2 \gamma}{\gamma + 1}} ((0,T) \times \Td; R^d)
		\]
		as $\ell \to \infty$, where 
\[
E[ \vr ] = \int_{X_D}  \vr [\hat{\vr}, \hat{\vu}, \hat{\mu} , \hat{\eta}, \hat{a}, \hat{\vc{g}}]\                \D \mathcal{L}_0	,\ 
E[ \vr \vu ] = \int_{X_D}  (\vr \vu) [\hat{\vr}, \hat{\vu}, \hat{\mu} , \hat{\eta}, \hat{a}, \hat{\vc{g}}]                \D \mathcal{L}_0;
\]
		
\item 
	\[
	\mathbb{E}_r \left[ \frac{1}{N} \sum_{n=1}^N \delta_{\vr^{h,n}} \right] \to 
	\mathbb{E}_r [\vr] \ \mbox{in} \ L^\gamma(\Td),\ 1 < r < \gamma, 
	\]
\[
\mathbb{E}_s \left[ \frac{1}{N} \sum_{n=1}^N \delta_{\vr^{h,n} \vu^{h,n}} \right] \to 
\mathbb{E}_s [\vr \vu] \ \mbox{in} \ L^{ \frac{2 \gamma}{\gamma + 1} }(\Td; R^d),\ 1 < s < \frac{2 \gamma}{\gamma + 1} 
\]	
as $\ell \to \infty$.	
		\end{itemize} 
	\end{Theorem}

\end{mdframed}

In accordance with the above definition, the empirical $r-$means are defined as minimizers:
\[
\mathbb{E}_r \left[ \frac{1}{N} \sum_{n=1}^N \delta_{\vr^{h,n}} \right] = {\rm arg} \min_{\rho \in L^\gamma (\Td) }
\frac{1}{N} \sum_{n=1}^N \left\| \vr^{h,n} - \rho \right\|^r_{L^\gamma (\Td)}, 
\]
and, similarly, 
\[
\mathbb{E}_s \left[ \frac{1}{N} \sum_{n=1}^N \delta_{\vr^{h,n} \vu_{h,n} } \right]  = {\rm arg} \min_{\vc{m} \in L^{\frac{2 \gamma}{\gamma + 1} } (\Td; R^d) }
\frac{1}{N} \sum_{n=1}^N \left\| \vr^{h,n} - \vm \right\|^s_{L^{\frac{2 \gamma}{\gamma + 1} } (\Td; R^d)}.
\]

Similarly results can be obtained if the $r$-barycenter is related to a closed convex subset of $X$.

\bigskip 

\centerline{Acknowledgement}

A part of this work was done during author's visit to the Basque Center for Applied Mathematics in Bilbao. 
The support and friendly atmosphere are acknowledged with thanks.

\def\cprime{$'$} \def\ocirc#1{\ifmmode\setbox0=\hbox{$#1$}\dimen0=\ht0
	\advance\dimen0 by1pt\rlap{\hbox to\wd0{\hss\raise\dimen0
			\hbox{\hskip.2em$\scriptscriptstyle\circ$}\hss}}#1\else {\accent"17 #1}\fi}



\end{document}